\documentclass[3p]{elsarticle}
\usepackage{slashbox}
\usepackage{amsmath,amsopn,amssymb,epsfig,amsthm}%,hyperref}
\usepackage{multirow}
\usepackage{color}

\definecolor{amethyst}{rgb}{0.6, 0.4, 0.8}
\definecolor{orange}{rgb}{1,0.5,0}

\newtheorem{Theorem}{Theorem}[section]
\newtheorem{Lemma}{Lemma}[section]

\newtheorem{example}{Example}[section]

\newproof{pot}{Proof of Theorem \ref{thm2}}
\newcommand{\bb}{\begin{bmatrix}}
\newcommand{\eb}{\end{bmatrix}}
\newcommand{\bl}[1]{\begin{list}{#1}{\usecounter{bean}}} \newcommand{\el}{\end{list}}
\newcommand{\bel}[1]{\begin{equation} \label{#1}} \newcommand{\eel}{\end{equation}}

\def\r2n2n{\mathbb{R}^{2n\times 2n}}
\def\c2n2n{\mathbb{C}^{2n\times 2n}}

\begin{document}

\date{}
\begin{frontmatter}
\title{On the maximal solution of the conjugate discrete-time algebraic Riccati equation}
%\tnoteref{t1}}%,t2}}
%

\author{Hung-Yuan Fan \corref{cor1}}%\fnref{fn1}}
\ead{hyfan@ntnu.edu.tw}
\address{Department of Mathematics, National Taiwan Normal University, Taipei 116325, Taiwan.}
\author{Chun-Yueh Chiang\corref{cor2}}%\fnref{fn2}}
\ead{chiang@nfu.edu.tw}
\address{Center for General Education, National Formosa
University, Huwei 632, Taiwan.}

\cortext[cor2]{Corresponding author}
%\fntext[fn1]{The first
%author was supported by the Ministry of Science and Technology of Taiwan under grant 110-2115-M-003-016.}
%\fntext[fn2]{ The second author was supported by the Ministry of Science and Technology of Taiwan under grant 109-2115-M-150-003MY2.}

\date{ }

\begin{abstract}
In this paper we consider a class of conjugate discrete-time Riccati equations,
arising originally from the linear quadratic regulation problem for discrete-time antilinear systems.
Under some mild assumptions and the framework of the fixed-point iteration,
a constructive proof is given for the existence of the maximal solution
to the conjugate discrete-time Riccati equation, in which the control
weighting matrix is nonsingular and its constant term is Hermitian.
Moreover, starting with a suitable initial matrix, we also show that the nonincreasing sequence generated by the fixed-point iteration
converges at least linearly to the maximal solution of the Riccati equation.
An example is given to demonstrate the correctness of our main theorem and provide considerable insights into the study of another meaningful solutions.
%It is to be expected that our theoretical results presented in this paper will
%play an important role in the oprimal control problems for discrete-time antilinear systems.
\end{abstract}


\end{frontmatter}

\section{Introduction} \label{sec1}
In this paper we are mainly concerned with the existence of the maximal solution to the conjugate
discrete-time algebraic Riccati equation (CDARE)
\begin{subequations}\label{cdare}
\begin{align}
  X = \mathcal{R}(X):= A^H \overline{X}A - A^H \overline{X}B(R+B^H \overline{X}B)^{-1} B^H \overline{X}A + H, \label{cdare-a}
     % &= A^H X(I + GX)^{-1}A + H, \label{dare-b}
\end{align}
{\color{black} or its equivalent expression}
 \begin{align}
 X = A^H \overline{X} (I+ G\overline{X})^{-1} A + H, \label{cdare-b}
\end{align}
\end{subequations}
where $A\in \mathbb{C}^{n\times n}$, $B\in \mathbb{C}^{n\times m}$, $R\in \mathbb{H}_m$ is nonsingular,
$H\in \mathbb{H}_n$, $I$ is the identity matrix of compatible size and
$G:= BR^{-1}B^H$, respectively. Here, $\mathbb{H}_\ell$ denotes the set of all $\ell \times \ell$ Hermitian matrices.
{ We use the notation $\widehat{M} := \overline{M} M$ throughout the paper, where $\overline{M}$ is the matrix obtained by taking the complex conjugate of each entry of a matrix $M$. The solution $X\in \mathbb{H}_n$ of the CDARE \eqref{cdare-a},
with $R_X:=R+B^H \overline{X}B$ being positive definite, is considered in this paper.}

For any $M,N\in \mathbb{H}_n$, the positive definite and positive semidefinite matrices are denoted by $M > 0$ and $M\geq 0$, respectively.
Moreover, we usually denote by $M \geq N$ (or $M \leq N$) if $M - N \geq 0$ (or $N - M \geq 0$) in the context.
{ For the sake of simplicity, the spectrum and spectral radius of $A\in \mathbb{C}^{n\times n}$ are
denoted by $\sigma (A)$ and $\rho (A)$, respectively. We say that the CDARE \eqref{cdare} has a maximal solution $X_M\in \mathbb{H}_n$
if $X_M\geq X$ {for any solution $X\in \mathbb{H}_n$ of the CDARE.} Thus, it follows from the definition of the maximality that $X_M$ is unique if it exists.} A matrix operator $f:\mathbb{H}_{n}\rightarrow\mathbb{H}_{n}$ is order preserving (resp. reversing) on $\mathbb{H}_{n}$ if $f(A)\geq f(B)$ (resp.  $f(A)\leq f(B)$) when $A\geq B$ and $A,B\in\mathbb{H}_{n}$.

%if $X =\mathcal{R}(X)$ and $X\in \mathbb{P}$,
%where $\mathcal{R}(\cdot)$ is the Riccati operator defined by \eqref{cdare-a}.

A class of CDAREs \eqref{cdare} arises from the linear quadratic regulation (LQR) optimal control problem for the discrete-time antilinear system
of the state space representation
\begin{equation} \label{antils}
  x_{k+1} = A\overline{x_k} + B\overline{u_k},\quad k\geq 0,
\end{equation}
where $x_k\in \mathbb{C}^n$ is the state response and $u_k\in \mathbb{C}^m$ is the control input.
 The main goal of this control problem is to find a state feedback control $u_k = -Fx_k$ such that the performance index
 \[ \mathcal{J}(u_k, x_0) := \sum_{k=0}^\infty (x_k^H H x_k + u_k^H R u_k)  \]
is minimized with $H\geq 0$ and $R > 0$. If the antilinear system \eqref{antils} is controllable,
it is shown in Theorem 12.7 of \cite{w.z17} that the optimal state feedback controller is $u_k^* = -R_{X_*}^{-1}B^H \overline{X_*}A x_k$ for $k\geq 0$ and thus
the minimum value of $\mathcal{J}(u_k^*, x_0) = x_0^H X_* x_0$ is achieved,
where $X_* \geq 0$ is the unique solution of the CDARE \eqref{cdare-a},
which is also called the {\em discrete-time algebraic anti-Riccati equation} \cite{Wu.16,w.z17}.

 Recently, some accelerated iterations have been proposed for solving the unique positive definite solution  of the CDARE \eqref{cdare} under positive definiteness assumptions with $G > 0$ and $H > 0$, see, e.g., \cite{l.c18,l.c20,m19} and the references therein. In addition, this numerical technique has also been utilized to some real-life applications recently,
see, e.g., \cite{l.w.g20,r.l.a20,r.l.e21}.
%
%With $H\in\mathbb{P}$, it is easy to prove that $-1\not\in\sigma(G\overline{H})$ and the CDARE~\eqref{cdare-a} can be transformed into
%a discrete-time algebraic Riccati equation (DARE) of the form
%\begin{equation}\label{dare2}
%X = \widetilde{\mathcal{R}} (X) =\widetilde{H}+\widetilde{A}^H  X \widetilde{A} - \widetilde{A}^H  X \widetilde{B}
%(\widetilde{R}+\widetilde{B}^H X \widetilde{B})^{-1} \widetilde{B}^H X\widetilde{A},
%\end{equation}
%where $\widetilde{\mathcal{R}} (X) := \mathcal{R}(\mathcal{R}(X)) $ and the coefficient matrices are given by
%\begin{align*}
%\widetilde{A}&= %\overline{A}A-\overline{A}B(R+B^H \overline{H} B)^{-1}B^H\overline{H}A
%=\overline{A}(I+G\overline{H})^{-1}A,\quad
%\widetilde{B} = \bb \overline{B} & \overline{A}B \eb,\\
%\widetilde{R} &= \overline{R}\oplus (R+B^H \overline{H} B),\quad
%%\widetilde{G} = \widehat{B} \widehat{R}^{-1} \widehat{B}^\ast=\overline{G}+\overline{A} (I+G\overline{H})^{-1}{G}\overline{A}^\ast,
%\widetilde{H} = H+{A}^H \overline{H}(I+G\overline{H})^{-1}{A}.
%\end{align*}
%Therefore, based on the semigroup property presented in \cite{l.c18,l.c20}, the unique solution $X>0$ to the CDARE \eqref{cdare-a}
%can be computed superlinearly from the standard DARE \eqref{dare2}.
%In addition, this numerical technique has also been utilized to some real-life applications recently,
%see, e.g., \cite{l.w.g20,r.l.a20,r.l.e21}.

In this work the concept of the maximal solution to discrete-time algebraic Riccati equations (DAREs) presented in Theorem~13.1.1 of \cite{l.r95} will be extended for the CDARE \eqref{cdare}
with $G,H\in \mathbb{H}_n$, and hopefully, it would play an important role in the LQR control problem for discrete-time antilinear systems.
The paper is organized as follows. In Section~\ref{sec2}, we introduce some useful notations and auxiliary lemmas
that will be used in our main result. In Section~\ref{sec3}, based on the framework of the fixed-point iteration (FPI),
the monotonicity property of a sequence generated by the FPI and the existence of the maximal solution to CDAREs \eqref{cdare}
will be addressed { in Theorem~\ref{thm3p1}} under reasonable hypotheses in \eqref{ma}. {We will give an example of the scalar CDARE \eqref{cdare} to illustrate the correctness of Theorem~\ref{thm3p1} in Section~\ref{sec4}.}
Finally, we conclude the paper in Section~\ref{conclusion}.

\section{Preliminaries}\label{sec2}

In this section we introduce some notations and auxiliary lemmas that will be used below. Firstly, let the conjugate Stein matrix operator $\mathcal{C}_A:\mathbb{H}_n\rightarrow\mathbb{H}_n$ associated with a matrix $A\in\mathbb{C}^{n\times n}$ be defined by
\begin{equation}\label{SME}
\mathcal{C}_A(X):= X-A^H \overline{X} A,
\end{equation}
for any $X\in\mathbb{H}_n$. In general, the operator $\mathcal{C}_A$ is neither order preserving nor order reversing. However, under the assumption that $\rho(\overline{A}A)<1$, its inverse operator $\mathcal{C}^{-1}_A$ is order preserving,
since $\mathcal{C}_A^{-1}(X)=\sum\limits_{k=0}^\infty ((\overline{A}A)^k)^H (X+A^H \overline{X} A) (\overline{A}A)^k \geq
\sum\limits_{k=0}^\infty ((\overline{A}A)^k)^H (Y+A^H \overline{Y} A) (\overline{A}A)^k=\mathcal{C}_A^{-1}(Y)$
for $X,Y\in \mathbb{H}_n$ with $X\geq Y$.

It will be clear later on that the matix operator $\mathcal{R}:\mbox{dom}(\mathcal{R})\rightarrow \mathbb{H}_n$ defined by
\eqref{cdare-a} plays an important role im Theorem~\ref{thm3p1} below, where $\mbox{dom}(\mathcal{R}):=\{X\in\mathbb{H}_n\, |\, \det(R_X)\neq 0\}$. Moreover, the following sets
\begin{equation}
 \mathcal{R}_\leq := \{X\in \mbox{dom}(\mathcal{R})\, |\, X  \leq  \mathcal{R}(X) \},\quad  \mathcal{R}_\geq := \{X\in \mbox{dom}(\mathcal{R})\, |\, X \geq \mathcal{R}(X) \}, \label{Rgeq} %\\
\end{equation}
will be used throughout the paper.
{Let $F_X = R_X^{-1}B^H \overline{X}A$ and $T_X = A-BF_X $ for any solution $X\in \mathbb{H}_n$ of the CDARE \eqref{cdare-a}.
Note that $T_X = (I+G\overline{X})^{-1}A$ and $\mathcal{R}(X) = A^H \overline{X} T_X + H= T_X^H \overline{X} A + H$, respectively.}

The following lemma characterizes { a useful identity} depending on the matrix operator $\mathcal{R}(\cdot)$ and
its associated conjugate Stein operator.
\begin{Lemma} \label{lem2p1} \color{black}
For any $X,Y\in\mbox{dom}(\mathcal{R})$, the following matrix identity
  \begin{align}\label{Req-b}
  X-\mathcal{R}(X)&=\mathcal{C}_{T_{Y}}(X)-H_{Y}+K(Y,X)
  \end{align}
  holds, where $K(Y,X) :=(F_Y-F_X)^H R_X(F_Y-F_X)$ and $H_{Y} :=H+K(Y,0)= H+F_{Y}^H R F_{Y}$.
\end{Lemma}
\begin{proof}
The proof is based on a direct calculation and the details are left to the reader.
%\begin{enumerate}
%  \item[(i)]
%  Observe that
% \[
% X-\mathcal{R}(X)=\Gamma_1(X)+\Gamma_2(X),
% \]
%where $\Gamma_1(X):=\mathcal{C}_A(X)-H$ and $\Gamma_2(X):=A^H \overline{X} B(R+B^H \overline{X} B)^{-1}B^H \overline{X} A$. A direct computation yields
%\begin{align*}
%&\Gamma_2(X)=F_X^H(R+B^H \overline{X} B)F_X=K_F(X)+[F^H(R+B^H \overline{X} B)F_X\\
%&+F_X^H(R+B^H \overline{X} B)F]-F^H(R+B^H \overline{X} B)F\\
%&=K_F(X)-F^H R F+[F^HB^H \overline{X} A+A^H \overline{X} BF-F^HB^H \overline{X} BF]\\
%&=K_F(X)-F^H R F+\mathcal{C}_{A_F}(X)-\mathcal{C}_{A}(X).
%\end{align*}
%We conclude that
%\begin{align*}
% X-\mathcal{R}(X)&=[\mathcal{C}_{A}(X)-H] + [K_F(X)-F^H R F+\mathcal{C}_{A_F}(X)-\mathcal{C}_{A}(X)] \\
% &=\mathcal{C}_{A_{F}}(X)-H_F+K_F(X).
%\end{align*}
%%
%\item[(ii)] The result is clearly true from \eqref{Req-a} with $F=F_{Y}$ and $A_{F_{Y}}\equiv T_{Y}$.
%\end{enumerate}
\end{proof}
%
%\begin{Lemma}\label{lem2p}
%Let $B\in \mathbb{C}^{n\times n}$ and $Q\geq 0$. If $X_0$ is a positive semidefinite solution of the conjugate Stein inequality $\mathcal{C}_B(X)\geq Q$,
%and $\mathrm{Ker}(Q)\cap\mathrm{Ker}(B^H\overline{Q}B)\subseteq \mathrm{Ker}(\overline{B}B-A)$ for some $A\in\mathbb{C}^{n\times n}$,
%then $\rho(\overline{B}B)\leq \max\{1,\rho(A)\}$. Furthermore, we have
%\begin{enumerate}
%  \item[(i)] $\rho(\overline{B}B)\leq 1$ if $\rho(A)\leq 1$.
%  %
%  \item[(ii)] $\rho(\overline{B}B)<1$ if $\rho(A)<1$ or $\mathrm{Ker}(Q)\cap E_\lambda(B)=\{0\}$ for some $\lambda\in\sigma(B)$.
%\end{enumerate}
%\end{Lemma}
%\begin{proof}
%Let $\lambda\in\sigma(\overline{B}B)$ and an eigenvector $x\in E_\lambda(\overline{B}B)$. Then
%$x^H \mathcal{S}_{\overline{B}B}(X_0) x=(1-|\lambda|^2)(x^H X_0 x) \geq x^H Q x +x^H B^H Q Bx\geq 0$. It is easily to see that
%\[
%E_\lambda(\overline{B}B)\subseteq\mathrm{Ker}(Q)\cap\mathrm{Ker}(B^H\overline{Q}B)
%\]
%for $\lambda\in\sigma(\overline{B}B)\cap\mathbb{D}^c$.
%
%If $x\not\in\mathrm{Ker}(Q)$, then $x\not\in\mathrm{Ker}(X_0)$ and thus $|\lambda|< 1$. Otherwise, $x\in \mathrm{Ker}(Q)$ implies that $Ax=Bx+(A-B)x=\lambda x$ and thus $|\lambda|\leq \rho(A)$. We complete the proof.
%
%\end{proof}
%

Let $\mathbb{T}:=\{X\in\mathbb{H}_n \,|\, \rho(\widehat{T}_X)<1\}$ and $\mathbb{P}:=\{X\in\mathbb{H}_n \,|\, R_X>0\}$.
{\color{black} A sufficient condition is provided in the following result to guarantee the stability property of a sequence generated by
the fixed-point iteration, which will be described in the proof of Theorem~\ref{thm3p1} later.}
   \begin{Lemma}\label{lem3p1}
Assume that $X_{\mathbb{T}}\in\mathbb{T}$.
If $X\in\mathbb{P}$ and $Y\geq 0$ both satisfy the matrix inequality
\begin{align}\label{3}
\mathcal{C}_{T_X}(Y)\geq K(X_{\mathbb{T}},X),
\end{align}
then $X\in\mathbb{T}$, i.e., $\rho (\widehat{T}_X) < 1$.
\end{Lemma}
\begin{proof}
The inequality \eqref{3} implies
\begin{align}\label{4}
\mathcal{S}_{\widehat{T}_X}(Y)\geq K(X_{\mathbb{T}},X)+T_X^H \overline{K(X_{\mathbb{T}},X)} T_X,
\end{align}
where $\mathcal{S}_{\widehat{T}_X}(Y) := Y - \widehat{T}_X^H Y \widehat{T}_X$. Assume that there exists a scalar $\lambda$ with $|\lambda|\geq 1$ such that $\widehat{T}_X x=\lambda x$ for some nonzero vector $x\in\mathbb{C}^n$.
Then we see that
\[
0\geq(1-|\lambda|^2)(x^H Y x)=x^H \mathcal{S}_{\widehat{T}_X}(Y) x \geq x^H K(X_{\mathbb{T}},X) x+x^H T_X^H \overline{K(X_{\mathbb{T}},X)} T_X x \geq 0.
\]
%For any $W\in\mathbb{H}_n$,
Note that $\mathrm{Ker}(K(X_{\mathbb{T}},X))=\mathrm{Ker}(F_{X_{\mathbb{T}}}-F_X)
\subseteq\mathrm{Ker}(T_{X_\mathbb{T}}-T_X)$ and
$\mathrm{Ker}( T_X^H\overline{K(X_{\mathbb{T}},X)} T_X) \subseteq\mathrm{Ker}(\overline{T_{X_\mathbb{T}}}T_X-\overline{T_X}T_X)$. This implies $T_X x = T_{X_{\mathbb{T}}}x,$ and
 $\overline{T_X}T_X x = \overline{T_{X_\mathbb{T}}}T_X x$.
Thus, $\widehat{T}_{X_\mathbb{T}}x=\overline{T_{X_\mathbb{T}}}T_{X_\mathbb{T}}x =\overline{T_{X_\mathbb{T}}}T_X x=\overline{T_X}T_X x=\widehat{T}_{X}x=\lambda x$,
which contradicts to the assumption, and hence we conclude that $\rho (\widehat{T}_X) < 1$.
\end{proof}

Let %$\mathcal{S}_{\geq} := \{ X\in \mathbb{H}_n\, |\, \mathcal{C}_{T_{X_\mathbb{T}}} (X)\geq H_{X_\mathbb{T}}\
%\mbox{for some } X_\mathbb{T}\in\mathbb{T}\}$
$\mathcal{S}_{\geq} :=\bigcup\limits_{X_\mathbb{T}\in\mathbb{T}} \{ X\in \mathbb{H}_n\, |\, \mathcal{C}_{T_{X_\mathbb{T}}} (X)\geq H_{X_\mathbb{T}}\}$
 if $\mathbb{T}\neq\emptyset$. For any $X_\mathbb{T}\in\mathbb{T}$, notice that the operator $\mathcal{C}_{T_{X_\mathbb{T}}}$ is a bijection
because $\rho (\widehat{T}_{X_\mathbb{T}}) < 1$ and its inverse operator $\mathcal{C}^{-1}_{T_{X_\mathbb{T}}}$ is order preserving. Therefore, $\mathcal{S}_\geq$ must be a nonempty subsets of $\mathbb{H}_n$ when $\mathbb{T}\neq\emptyset$.

{\color{black} In order to deduce our main theorem presented in the next section, we assume that
\begin{equation}\label{ma}
 \mathbb{T}\neq\emptyset\quad \mbox{and}\quad \mathcal{R}_{\leq}\cap\mathbb{P}\neq\emptyset
\end{equation}
in the following content of the paper.
Starting with $X_0\in \mathcal{S}_\geq$,}
the monotonicity property of a sequence $\{X_k\}_{k=0}^\infty$ generated by the following FPI of the form
\begin{equation} \label{fpi}
X_{k+1} = \mathcal{R}(X_k), \quad k\geq 0,
\end{equation}
will be ensured by the following lemma.
\begin{Lemma}\label{lem3p2}
\par\noindent
%Assume that the assumptions in \eqref{ma} hold.
{\color{black} If $X_\star\in\mathcal{S}_{\geq}$, then the following statements hold:}
\begin{enumerate}
  \item[(i)]
 $X_\star\geq X_\mathbb{P}$ for any $X_\mathbb{P}\in \mathcal{R}_\leq\cap\mathbb{P}$.
 \item[(ii)]
 $X_\star-\mathcal{R}(X_\star)\geq K(X_{\mathbb{T}},X_\star)$ for some $X_{\mathbb{T}}\in\mathbb{T}$. Furthermore, we have  $\mathcal{S}_{\geq}\subseteq \mathcal{R}_\geq\cap\mathbb{P}$.
 \end{enumerate}
\end{Lemma}
\begin{proof}
Since $X_\star\in \mathcal{S}_{\geq}$, there is a matrix
$X_\mathbb{T} \in \mathbb{T}$ associated with $X_\star$ such that $\mathcal{C}_{T_{X_\mathbb{T}}} (X_\star)\geq H_{X_\mathbb{T}}$.
 \begin{enumerate}
\item[(i)] {\color{black}
For any $X_\mathbb{P}\in \mathcal{R}_\leq\cap\mathbb{P}$, it follows from \eqref{Req-b} that
\[
0\geq X_\mathbb{P}-\mathcal{R}(X_\mathbb{P})=\mathcal{C}_{T_{X_\mathbb{T}}}(X_\mathbb{P})-H_{X_\mathbb{T}}+K(X_\mathbb{T},X_\mathbb{P}).
\]
Since $X_{\mathbb{P}}\in \mathbb{P}$, $K(X_\mathbb{T},X_\mathbb{P})\geq 0$ and hence we deduce that
%\begin{align}\label{leq1}
$H_{X_\mathbb{T}}\geq\mathcal{C}_{T_{X_\mathbb{T}}}(X_\mathbb{P})$.
%\end{align}
Therefore,
\[
    \mathcal{C}_{T_{X_\mathbb{T}}}(X_\star-X_\mathbb{P})
    =\mathcal{C}_{T_{X_\mathbb{T}}}(X_\star)-\mathcal{C}_{T_{X_\mathbb{T}}}(X_\mathbb{P})
    \geq H_{X_\mathbb{T}}-\mathcal{C}_{T_{X_\mathbb{T}}}(X_\mathbb{P})\geq 0
\]
and then  $X_\star\geq X_\mathbb{P}$. Furthermore, we also show that $\mathcal{S}_{\geq}\subseteq \mathbb{P}$,
          since $R_{X_{\star}}\geq R_{X_{\mathbb{P}}}  > 0$.

\item[(ii)]
From \eqref{Req-b}, note that
$X_\star-\mathcal{R}(X_\star)=\mathcal{C}_{T_{X_\mathbb{T}}}(X_\star)-H_{X_\mathbb{T}}+K(X_\mathbb{T},X_\star)\geq K(X_\mathbb{T},X_\star)\geq 0$
because $X_\star\in \mathbb{P}$, and we thus conclude that $\mathcal{S}_{\geq}\subseteq  \mathcal{R}_\geq \cap \mathbb{P}$.}
%Again from \eqref{Req-b} one obtains that
%\[
%X_\star-\mathcal{R}(X_\star)=\mathcal{C}_{T_{X_\mathbb{T}}}(X_\star)-H_{X_\mathbb{T}}+K(X_\mathbb{T},X_\star)\geq K(X_\mathbb{T},X_\star).
%\]
\end{enumerate}
\end{proof}

\section{The maximal solution of the CDARE} \label{sec3}

In this section, the existence of the maximal solution to the CDARE \eqref{cdare}
will be established iteratively, utilizing the framework of the FPI \eqref{fpi}
with an appropriate initial matrix $X_0\in \mathbb{H}_n$, {\color{black} under the assumptions in \eqref{ma}. The main results of this paper are
summarized in the following theorem.}
 %Let sequence $\{X_k\}_{k=0}^\infty$ be
%generated by the FPI \eqref{fpi} with $X_0\in\mathcal{S}_\geq$.
\begin{Theorem} \label{thm3p1}
{\color{black} If the assumptions in \eqref{ma} are fulfilled, then the maximal solution $X_M$ of the CDARE \eqref{cdare} exists.}
Furthermore, the following statements hold:
\begin{enumerate}
  \item[(i)]
The sequence $\{X_k\}_{k=0}^\infty$ generated by the FPI \eqref{fpi} with $X_0\in\mathcal{S}_\geq$ is well-defined. Moreover,
   $X_k\in\mathcal{S}_\geq \cap \mathbb{T}\subseteq \mathcal{R}_\geq \cap \mathbb{P} \cap \mathbb{T}$ for all $k\geq 0$.
  \item[(ii)] $X_k\geq X_{k+1}\geq X_\mathbb{P}$ for all $k\geq 0$ and $X_\mathbb{P}\in \mathcal{R}_{\leq} \cap \mathbb{P}$.
  \item[(iii)] The sequence $\{X_k\}_{k=0}^\infty$  converges at least linearly to $X_M$, which is the maximal element of the set $\mathcal{R}_\leq\cap \mathbb{P}$
and satisfies $\rho (\widehat{T}_{X_M}) \leq 1$, with the rate of convergence
  \[  \limsup_{k\rightarrow \infty} \sqrt[k]{\|X_k - X_{M}\|}\leq \rho (\widehat{T}_{X_M}) \]
 whenever $X_M\in\mathbb{T}$.
%
%\item[(iii)] For each $k\geq 0$, $X_k\in\mathbb{T}\cap\mathbb{P}$.
\end{enumerate}
\end{Theorem}
\begin{proof}
Let $X_{\mathbb{P}} \in \mathcal{R}_\leq \cap \mathbb{P}$.
The existence of $X_M$ can be constructed by the following discussion.
\begin{enumerate}
  \item[(i)]
   {\color{black} This result will be proven by induction.} Since $X_0 \in \mathcal{S}_\geq$, it follows from Lemma \ref{lem3p2} that $X_0 \geq X_{\mathbb{P}}$ and $X_0 \in \mathcal{R}_\geq \cap \mathbb{P}$. Thus, $X_1 = \mathcal{R} (X_0)$ is well-defined with $X_0 \geq X_1$. Furthermore, we also have
   \begin{align*}
 \mathcal{C}_{T_{X_0}}(X_0-X_\mathbb{P}) &= \mathcal{C}_{T_{X_0}}(X_0)
 -\mathcal{C}_{T_{X_0}}(X_\mathbb{P})= H_{X_0}+X_0-\mathcal{R}(X_0)-\mathcal{C}_{T_{X_0}}(X_\mathbb{P})\\
 &\geq H_{X_0}+K(X_\mathbb{T},X_0)-(H_{X_0}+X_\mathbb{P}-\mathcal{R}(X_\mathbb{P})-K(X_0, X_\mathbb{T}))\\
 &= K(X_\mathbb{T},X_0)-(X_\mathbb{P}-\mathcal{R}(X_\mathbb{P}))+K(X_0,X_\mathbb{P})
 \geq K(X_\mathbb{T},X_0),
\end{align*}
for some $X_{\mathbb{T}} \in \mathbb{T}$. Then $X_0\in \mathbb{T}$ follows immediately from Lemma \ref{lem3p1} since
  $\rho(\widehat{T}_{X_{\mathbb{T}}})<1$. Assume that $X_k \in \mathcal{S}_\geq \cap \mathbb{T}$ for some positive integer $k$.
  Again, it follows from Lemma \ref{lem3p2} that $X_k \geq X_{\mathbb{P}}$, $X_k \in \mathcal{R}_\geq \cap \mathbb{P}$ and
 $X_{k+1} = \mathcal{R} (X_k)$ is well-defined with $X_k \geq X_{k+1}$.
  In addition, from Lemma~\ref{lem2p1} we have
\begin{subequations}
\begin{align}%\label{eq}
  X_k-\mathcal{R}(X_k)&=\mathcal{C}_{T_{X_k}}(X_k)-H_{X_k},\label{eq1}\\
  X_{k+1}-\mathcal{R}(X_{k+1})&=\mathcal{C}_{T_{X_{k+1}}}(X_{k+1})-H_{X_{k+1}},\label{eq2}\\
  X_{k+1}-\mathcal{R}(X_{k+1})&=\Delta_k+K(X_{k},X_{k+1}),\label{eq3}\\
  X_\mathbb{P}-\mathcal{R}(X_\mathbb{P}) &=(\mathcal{C}_{T_{X_{k+1}}}(X_\mathbb{P})-H_{X_{k+1}})+K(X_{k+1},X_\mathbb{P}), \label{eq4}
\end{align}
\end{subequations}
  where $\Delta_k := \mathcal{C}_{T_{X_k}}(X_{k+1})-H_{X_k}$.
  From \eqref{eq1} we see that
\begin{align*}
\Delta_k %&:= \mathcal{C}_{T_{X_k}}(X_{k+1})-H_{X_k} \\
&=\mathcal{C}_{T_{X_k}}(X_{k+1})+X_k-X_{k+1}-\mathcal{C}_{T_{X_k}}(X_{k})\\
&=-\mathcal{C}_{T_{X_k}}(X_k-X_{k+1})+X_k-X_{k+1}=T_k^H(\overline{X_k-X_{k+1}})T_k\geq 0.
\end{align*}
  Thus, it implies that $X_{k+1}\in \mathcal{S}_\geq$ because $X_k\in \mathbb{T}$, and thus $X_{k+1} \geq X_{\mathbb{P}}$ follows from Lemma \ref{lem3p2}.
  Moreover, from \eqref{eq2}, \eqref{eq3} and \eqref{eq4}, the positive semidefinite matrix $Y:=X_{k+1} - X_{\mathbb{P}}$ satisfies the conjugate Stein matrix inequality
  \begin{align*}
 \mathcal{C}_{T_{X_{k+1}}}(Y) &= \mathcal{C}_{T_{X_{k+1}}}(X_{k+1})-
 \mathcal{C}_{T_{X_{k+1}}}(X_\mathbb{P})\\
 &=\left(X_{k+1}-\mathcal{R}(X_{k+1})+H_{X_{k+1}}\right)-\left((X_\mathbb{P}-\mathcal{R}(X_\mathbb{P}))+H_{X_{k+1}}
 -K(X_{k+1},X_\mathbb{P}) \right)\\
 &=\Delta_k+K(X_{k},X_{k+1})+H_{X_{k+1}}-(X_\mathbb{P}-\mathcal{R}(X_\mathbb{P}))-H_{X_{k+1}}
 +K(X_{k+1},X_\mathbb{P})\\
 &=\Delta_k+K(X_{k},X_{k+1})-(X_\mathbb{P}-\mathcal{R}(X_\mathbb{P}))+K(X_{k+1},X_\mathbb{P}) \geq K(X_{k},X_{k+1}),
\end{align*}
and hence $X_{k+1} \in \mathbb{T}$ {\color{black} follows from Lemma \ref{lem3p1} and $X_k \in \mathbb{T}$. That is,
we have deduced that $X_{k+1} \in \mathcal{S}_\geq \cap \mathbb{T}$ and this completes the inductive proof.}
  \item[(ii)] The statement follows immediately from the proof of (i) and Lemma \ref{lem3p2}.
  \item[(iii)] Since it follows from the statement (ii) that $\{X_k\}_{k=0}^\infty$ is a nonincreasing sequence of Hermitian matrices and bounded below by all elements
  of the set $\mathcal{R}_\leq\cap \mathbb{P}$, $X_M := \lim\limits_{k \rightarrow \infty} X_k$ must be the maximal element of $\mathcal{R}_\leq\cap \mathbb{P}$
  such that $X_M = \mathcal{R} (X_M)$ and $\rho (\widehat{T}_{X_M}) \leq 1$.
  Moreover, the proof of the rate of convergence for $\{X_k\}_{k=0}^\infty$ follows from the Appendix A.1 of \cite{l.c18}.
\end{enumerate}
\end{proof}

\section{An example for the CDARE} \label{sec4}
In this section we shall give an example to illustrate the correctness of Theorem~\ref{thm3p1} presented
in the previous section. %In our numerical experiments shown below, we will measure the normalized residual
%\[  NRes(Z) : = \frac{\|Z - \mathcal{R}(Z) \|}
%{\|Z\| + \|A^H \overline{Z} A\| + \|A^H \overline{Z} B R_Z^{-1} B^H \overline{Z} A\| + \|H\|}  \]
%for each quantity $Z$ computed by the FPI \eqref{fpi}, where $\|\cdot \|$ denotes the matrix $2$-norm.
% In our numerical results, the FPI \eqref{fpi} was terminated when $NRes \leq 1.0\times 10^{-15}$.
%
%All numerical experiments were
%performed on ASUS laptop (ROG GL502VS-0111E7700HQ), using Microsoft Win-
%dows 10 operating system and MATLAB Version R2019b, with Intel Core i7-7700HQ
%CPU and 32 GB RAM.

\begin{example} \label{ex1}   \em
We consider a scalar CDARE \eqref{cdare} of the form
\begin{equation} \color{black}
  x  = |a|^2 \bar{x} -\frac{|a|^2 \bar{x}^2 |b|^2}{r + |b|^2\bar{x}} + h %\nonumber \\
   = \frac{|a|^2 \bar{x}}{1+g\bar{x}} + h, \label{scdare}
\end{equation}
where $a,b\in \mathbb{C}$, $r,h\in \mathbb{R}$ with $r+|b|^2\bar{x} > 0$ and $g := |b|^2 / r$ with $r\neq 0$. Without loss of generality, we assume that
$\color{black}  |a| > 0$, $r>0$ and thus $g>0$.
From \eqref{scdare} and $\color{black} 1+g\bar{x} > 0$, we obtain
\[ g|x|^2 + x - (|a|^2 + gh)\bar{x} - h = 0,  \]
which has two solutions $\color{black} x_M, x_m \in \mathbb{H}_1 = \mathbb{R}$ satisfying
\[  x_M := \frac{-(1 - |a|^2 - gh) + \sqrt{D}}{2g}\quad \mbox{and}\quad  x_m := \frac{-(1 - |a|^2 - gh) - \sqrt{D}}{2g}\]
if %$1+|a|^2+gh \neq 0$ and
the discriminant
$D := (1 - |a|^2 - gh)^2 + 4gh \geq 0$. Note that
$D\geq 0$ if and only if $h\geq h_M:=\frac{-(1-|a|)^2}{g}$ or $h\leq h_m:=\frac{-(1+|a||)^2}{g}$.

{\color{black} Let $\hat{t}_x : = \widehat{T}_x = \frac{|a|^2}{(1+gx)^2} $ for $x=x_M$ or $x = x_m$. Some facts are listed without proof as follows:}
\begin{enumerate}
  \item $\mathbb{T}=(-\infty,\frac{-|a|-1}{g})
      \cup(\frac{|a|-1}{g},\infty)\neq\emptyset$ %if $\color{black} g > 0$
      and $\mathcal{R}_{\leq}\cap\mathbb{P}=[x_m,x_M]\cap(\frac{-1}{g},\infty)\neq\emptyset$ if $h\geq h_M$. Moreover, $1+gx_M\geq1+gx_m>0$ and  %$\mathcal{S}_\geq\subseteq (\frac{-1}{g},x_m]\cup[x_M,\infty)$ if $h\geq h_M$.
      $\color{black} \mathcal{S}_\geq\subseteq [x_M,\infty)$ if $h\geq h_M$.
      When $h\leq h_m$, we have $1+gx_m<1+gx_M<0$ and $\mathcal{R}_{\leq}\cap\mathbb{P}=\emptyset$.
  \item $\hat{t}_{x_M}\hat{t}_{x_m}=1$. Furthermore,
   $\hat{t}_{x_M}< 1 < \hat{t}_{x_m}$ if $h> h_M$,
    $\hat{t}_{x_M}=1=\hat{t}_{x_m}$ if $h=h_M$ or $h=h_m$, and
   $\hat{t}_{x_M}> 1 > \hat{t}_{x_m}$ if $h< h_m$.
  \item {\color{black} For $x_0\in\mathcal{S}_{\geq}\backslash \{x_M\}$, the equivalent expressions of the sequence $\{x_k\}_{k=0}^\infty$ generated by the FPI \eqref{fpi}  can be rewritten in the form}
  \begin{subequations}\label{exp}
\begin{align}
x_{k}
&=x_M+\dfrac{x_M-x_m}{s\hat{t}_{x_M}^{-k}-1}\ \mbox{ if}\  0 < \hat{t}_{x_M} < 1, \label{exp1}\\
&=x_M+\dfrac{|a|(x_0-x_M)}{g(x_0-x_M)k+|a|}\ \mbox{ if}\ \hat{t}_{x_M}=1, \label{exp2}
 \end{align}
  \end{subequations}
for all $k\geq 0$, where $s:=\frac{x_0-x_m}{x_0-x_M}$.
\end{enumerate}

{\color{black} Based on these facts, the convergence behavior of the sequence presented in \eqref{exp} is summarized below for different cases depending on $h$.}
\begin{itemize}
  \item [(i)] {\color{black} If $h>h_M$ and $0 < \hat{t}_{x_M} < 1$, then the assumptions in \eqref{ma} hold and the sequence $\{x_k\}_{k=0}^\infty$ in the expression~\eqref{exp1} converges
      decreasingly to $x_M$, which is the maximal solution of the CDARE \eqref{scdare}.} Furthermore, the convergence is linearly since $\lim\limits_{k\rightarrow\infty}\sqrt[k]{|x_k-x_M|}=\hat{t}_M<1$. {\color{black} When $x_0\in(x_m,x_M)$,
      we see that $s<0$ and thus $\{x_k\}_{k=0}^\infty$ converges increasingly to $x_M$ with linear convergence even though $x_0\in\mathcal{S}_{\geq}$ is
      not true in this case.}

  \item [(ii)] {\color{black} If $h=h_M$, then $x_M = x_m$ with $\hat{t}_{x_M} = 1$ and the assumptions in \eqref{ma} hold.
  In this case we see that the sequence $\{x_k\}_{k=0}^\infty$ in the expression \eqref{exp2} converges decreasingly to $x_M$ with
      sublinear convergence, since $\lim\limits_{k\rightarrow\infty}\frac{|x_{k+1}-x_M|}{|x_k-x_M|}=\hat{t}_M=1$.}
  \item [(iii)] {\color{black} If $h<h_m$, then the assumptions in \eqref{ma} do not hold and $\hat{t}_M>1$.}
  However, from \eqref{exp1} it can be shown that $x_k\rightarrow x_M+\frac{x_M-x_m}{-1}=x_m$ as $k\rightarrow\infty$. Namely,
       the sequence $\{x_k\}_{k=0}^\infty$ converges linearly to the minimal solution {\color{black} of the CDARE \eqref{scdare}.
       Further investigations are left as a topic of our future work for this case.}
\end{itemize}
\end{example}
\section{Concluding remarks} \label{conclusion}

In this paper we mainly deal with a class of conjugate discrete-time Riccati equations,
arising originally from the LQR control problem for discrete-time antilinear systems. In this case, the design of the optimal controller usually depends on
the existence of a unique positive semidefinite optimizing solution of CDAREs \eqref{cdare-a} with $R>0$ and $H\geq 0$,
if the antilinear system is assumed to be controllable.

Analogous to Theorem 13.1.1 of \cite{l.r95} for standard DAREs, based on the framework of the fixed-point iteration, we have proved the existence of
the maximal solution to the CDARE \eqref{cdare-a}, with nonsingular $R$ and $H\in  \mathbb{H}_n$, under some weaker assumptions defined as in \eqref{ma}. Moreover, starting from $X_0 \in \mathcal{S}_\geq$, the monotonicity of the sequence $\{X_k\}_{k=0}^\infty$
generated by the FPI \eqref{fpi} and the stability of its corresponding sequence
$\{\widehat{T}_{X_k} \}_{k=0}^\infty$ have also been deduced under the same assumptions, respectively.
We believe that our theoretical results would be useful in the LQR control problem,
or even the state-feedback stabilization problem, for discrete-time antilinear systems.

It is only to be expected that the existence of the minimal solution or other extremal solutions of the CDARE \eqref{cdare}
will be investigated in the next work, and
{ it also leads to our future work that how to apply the accelerated techniques presented in \cite{l.c18,l.c20}
for solving the extremal solutions to the CDARE \eqref{cdare}. }
\section*{Acknowledgment}
This research work is partially supported by the Ministry of Science and Technology and the National Center for Theoretical Sciences in Taiwan. The first author (Hung-Yuan Fan) would like to thank the support from the Ministry of Science and Technology of Taiwan
under grants MOST 110-2115-M-003-016, and the corresponding author (Chun-Yueh Chiang) would like to thank the support from the Ministry of Science and Technology of Taiwan under the grant MOST 109-2115-M-150-003-MY2.
%\bibliographystyle{abbrv}
%\bibliography{cdare}

\end{document}